\documentclass[a4paper,12pt]{article}
\usepackage{amsmath,amsfonts,amssymb,amscd}
\usepackage{hyperref}
\input{xypic}

\usepackage{multicol}

\def\s{\scriptstyle }
\def\a{\alpha}
\def\l{\lambda}

\def\x{\mathbf x}
\def\y{\mathbf y}

\def\mfS{{\mathfrak S}}

\def\Tab{\mathfrak{ T\hspace{-0.07 em}a\hspace{-0.05 em}b}}

\def\C{{\mathbb C}}

\def\moins{\raise 1pt\hbox{{$\scriptstyle -$}}}
\def\plus{\raise 1pt\hbox{{$\scriptstyle +$}} }

\newtheorem{theorem}{Theorem}
\newtheorem{proposition}[theorem]{Proposition}

\newdimen\carresize
\newdimen\thickness
\def\CARRE#1{\hbox{\vrule width \thickness
   \vbox to \carresize{\hrule height \thickness\vss
      \hbox to \carresize{\hss#1\hss}
   \vss\hrule height\thickness}
\unskip\vrule width \thickness}
\kern-\thickness}
\def\vsquare#1{\vbox{\CARRE{$#1$}}\kern-\thickness}

\def\smallyoung#1{%
  \carresize=12pt%
  \thickness=0.5pt%
  \vcenter{%
    \vbox{\smallskip\offinterlineskip%
      \halign{&\vsquare{##}\cr #1}}}}

\def\young#1{%
  \carresize=16pt%
  \thickness=0.5pt%
  \vcenter{%
    \vbox{\smallskip\offinterlineskip%
      \halign{&\vsquare{##}\cr #1}}}}

\newdimen\unit
\def\o{$\scriptscriptstyle{{\rm o}}$}
\def\grape(#1,#2)#3{\raise#2\unit\rlap{\kern#1\unit #3}\ignorespaces}

\def\gg{{\unit=1mm
\hbox {\grape(2,1.2){'}
       \grape(1,2)\o
       \grape(2,2)\o
       \grape(3,2)\o
       \grape(1.5,1)\o
       \grape(2.5,1)\o
       \grape(2,0)\o\
}\kern 3.5 \unit}}

\def\gfill{\leaders\hbox to 1.2em{\hss\gg\hss}\hfill}

\def\frise{\centerline{\hbox to 8cm{\gfill} }\bigskip}

\def\Pol{\mathfrak{P\hspace{-0.07 em}o\hspace{-0.07 em}l}}
\def\tc{\clubsuit}   
\def\tr{\spadesuit}  
\def\c{\mathfrak c}

\newdimen\vcadre\vcadre=0.2cm 
\newdimen\hcadre\hcadre=0.2cm 
\def\GrTeXBox#1{\vbox{\vskip\vcadre\hbox{\hskip\hcadre%
      $#1$%
   \hskip\hcadre}\vskip\vcadre}}
\def\arx#1[#2]{\ifcase#1 \relax \or%
  \ar @{-}[#2]  \or%
  \ar @2{-}[#2] \or%
  \ar @{--}[#2] \or%
  \ar @2{.}[#2] \or%
  \ar @2{~}[#2]  \fi}
\begin{document}

\date{}

\title{\bf Idempotents with polynomial coefficients}
\author{Alain Lascoux}

\maketitle

\begin{abstract}
We combine Young idempotents in the group algebra of the symmetric group
with the action of the symmetric group on 
 products of Vandermonde determinants 
to obtain idempotents with polynomial coefficients.
\end{abstract}

\vspace*{-10pt}
\hfill \emph{ Key words}:  Young idempotents, symmetric group

\medskip
Young described a complete set of idempotents in the group
algebra of the symmetric group. However, given an explicit idempotent,
how does one connect it with Young's work? 

In this text, I take the example of (quasi)-idempotents of the type
$$ \sum_{\sigma\in\mfS_n} (x_{\sigma_1} -x_{\sigma_2}) \, \sigma $$
in the group algebra of the symmetric group $\mfS_n$
with polynomial coefficients.

Given a list of integers $\alpha= [\a_1,\dots,\a_m]$, let 
$$ \Delta(\alpha) = \prod_{1\leq i\le j\leq m} (x_{\a_i}-x_{\a_j}) \, .$$
Given a partition $\l=[\l_1,\dots,\l_r]$ of 
weight $n=|\l|= \l_1\plus \dots \plus \l_r$, and a permutation 
$\sigma\in\mfS_n$, define
$$ \Delta(\sigma)= \Delta(\sigma_1,\dots,\sigma_{\l_1})  
    \Delta(\sigma_{\l_1+1},\dots, \sigma_{\l_1+\l_2}) \dots
   \Delta(\sigma_{\l_1+\dots+\l_{r-1}+1},\dots, 
           \sigma_{\l_1+\dots+\l_r}) \, ,$$      
which corresponds to cutting $\sigma$ into blocks of
respective lengths $\l_1,\l_2,\dots$.       

We shall consider the element
$$ \Omega_\l := \sum_{\sigma\in\mfS_n} \Delta_\l(\sigma)\, \sigma
  \ \in \ \Pol(x_1,\dots, x_n)\big[ \mfS_n \big] \, ,$$
and show that this element is quasi-idempotent,
i.e. that there exists a (non-zero) polynomial $c_\l$ such that
$\Omega_\l \Omega_\l = c_\l \Omega_\l$.   

Let us evoke in a few words the work of Young. 
Young described a linear basis $\{t_{tu}\}$ of 
the group algebra $\C[\mfS_n]$, indexed by pairs of standard 
tableaux  $t,u$ of the same shape of $n$ boxes.   
These $n!$ elements satisfy
$$  e_{t,u} e_{v,w} = \delta_{u,v}\, e_{t,w} \, .$$
Young basis can be characterized as the family of elements in $\C[\mfS_n]$
which are left and right eigenfunctions of the Jucys-Murphy elements
\cite{OkounkovVershik, LivreRep}.

The central idempotent corresponding to a partition $\l$ of $n$
is $e_\l= \sum_t e_{t,t}$, sum over all standard tableaux
of shape $\l$. 
In particular, when $\l$ has only one part, the idempotent
$e_n= (n!)^{-1} \sum \sigma$ is called the \emph{trivial idempotent}.

The unit of $\C[\mfS_n]$ decomposes as 
$ \sum_{\l:\, |\l|=n} e_\l =\sum_t e_{t,t}$, 
sum over all standard tableaux.  For example, for $n=2$,
there are two tableaux and the corresponding decomposition is
$1= \frac{1}{\s 2} (1\plus s_1) + \frac{1}{\s 2} (1\moins s_1)$.
We more generally write $s_1,\dots, s_{n-1}$ for
the \emph{simple transpositions},
which  generate $\mfS_n$. More generally, given a partition $\l$,
the group  $\mfS_\l$ generated by all $s_i$, $i\not\in 
\{ \l_1, \l_1\plus \l_2, \l_1\plus \l_2\plus \l_3,\dots \}$
is called a \emph{Young subgroup}.

Given a partition $\l$, let $\Tab(\l)$ be the set of standard
tableaux of columns of respective lengths $\l_1,\l_2,\ldots$.
The \emph{first tableau}  $\tc$ 
is the tableau filled by consecutive integers
in each column, and similarly, for the \emph{last tableau} $\tr$
one requires consecutive integers in each row.

Idempotents $e_{t,t}$ corresponding to the same shape
are obtained from any of them 
by conjugation by factors of the type $(s_i+c)$.

For example, for the shape $[2,2,1]$, the five idempotents are 
obtained recursively from the first idempotent $e_{\tc,\tc}$,
according to the following graph.

$\hspace*{-15pt} 
 \begin{array}{c}
 \xymatrix@R=0.5cm@C=-2cm{
 &  &  & *{\GrTeXBox{ \tc = \smallyoung{3\cr 2&5\cr 1&4\cr}}}
\arx5[d]& \\
 &  &  & *{\GrTeXBox{ \begin{array}{c}
 \smallyoung{4\cr 2&5\cr 1&3\cr} \\ \noalign{\kern 3pt}
      (s_3\moins \frac{1}{3})
     \frac{e_{\tc,\tc}}{1\moins \frac{1}{9}}
      (s_3\moins \frac{1}{3})  \end{array}
}}\arx1[ld]\arx3[rd]& \\
 &  & *{\GrTeXBox{
 \begin{array}{c} \smallyoung{4\cr 3&5\cr 1&2\cr} \\ \noalign{\kern 3pt}
 (s_2\moins \frac{1}{2})(s_3\moins \frac{1}{3})
   \frac{e_{\tc,\tc}}{(1\moins \frac{1}{9})(1\moins \frac{1}{4})}
(s_3\moins \frac{1}{3}) (s_2\moins \frac{1}{2})
  \end{array}
   }}\arx3[rd]&  & *{\GrTeXBox{ \begin{array}{c}
 \smallyoung{5\cr 2&4\cr 1&3\cr} \\  \noalign{\kern 3pt}
   (s_4\moins \frac{1}{2}) (s_3\moins \frac{1}{3})
  \frac{e_{\tc,\tc}}{(1\moins \frac{1}{9})(1\moins \frac{1}{4})}
(s_3\moins \frac{1}{3})(s_4\moins \frac{1}{2})
   \end{array}
}}\arx1[ld]& \\
 &  &  & *{\GrTeXBox{
  \tr =\smallyoung{5\cr 3&4\cr 1&2\cr}
}}& \\
} \\
e_{\tr,\tr} = (s_2\moins \frac{1}{2})
(s_4\moins \frac{1}{2}) (s_3\moins \frac{1}{3})
  \frac{e_{\tc,\tc}}{(1\moins \frac{1}{9})(1\moins \frac{1}{4})
     (1\moins \frac{1}{4})   }
(s_3\moins \frac{1}{3})(s_4\moins \frac{1}{2})
(s_2\moins \frac{1}{2})
\end{array}
$

\bigskip
The general rule to write the factors $(s_i\plus c)$ translates
in terms of group algebra Young's matrices of representation
\cite{LivreRep}.   Define the \emph{content} $c(i)$ of the integer $i$
located in the $j$-th column and $k$-th row of a tableau $t$
to be $j\moins k$. Then, given any $i$ such that the image
$s_i t$  of $t$ under the 
exchange of $i$ and $i\plus 1$ is still a tableau,
 one has, for any $u$ of the same shape,
\begin{eqnarray*}
e_{s_it,u} &=&  \frac{1}{\sqrt{1-\frac{1}{(c(i) -c(i+1))^2   } }}
\left( s_i +\frac{1}{c(i) -c(i\plus 1)}  \right)\, 
    e_{t,u}   \\
 e_{u,s_it}  &=&  e_{u,t}
  \left( s_i +\frac{1}{c(i) -c(i\plus 1)}  \right)
 \frac{1}{\sqrt{1-\frac{1}{(c(i) -c(i+1))^2   } }}
\end{eqnarray*}

For a pair of tableaux $t,u$ in $\Tab(\l)$, let 
$s_is_j\cdots s_k$ be a product (it does  to need to be of minimum length)
such that $s_i s_j\cdots s_k\, u =t$, denote by 
$\sigma_{t,u}$ this product, and by $\zeta(t,u)$ the corresponding
product of factors 
$\frac{1}{\sqrt{1-\frac{1}{(c(i) -c(i+1))^2   } }} 
\left( s_i +\frac{1}{c(i) -c(i\plus 1)}  \right)$.   
Beware that the contents are relative to variable tableaux in a sequence.
The above 
 products are better understood when using Yang-Baxter graphs
and spectral vectors \cite{SgPfaff}.

Given $\l$, let $Y_{\tc}= \Delta_\l(1,\ldots, n)$. The module
$V_\l$ generated by the action of $\mfS_n^x$ on $Y_{\tc}$ is
irreducible and called \emph{Specht module}\cite{Fulton}.
The \emph{Specht basis} is $\{ \sigma_{t,\tc}^x Y_{\tc}, \, 
t\in\Tab(\l) \}$,
while the Young basis is
$$  \{ Y_t:= e^x_{t,\tc} Y_{\tc} = 
    \zeta(t,\tc) Y_{\tc}, \, t\in\Tab(\l) \} \, .$$

\medskip
Given two copies $\mfS_n^1,\mfS_n^2$ of
$\mfS_n$, one defines the \emph{diagonal group} $\mfS_n^{12}$ to
be the other copy of $\mfS_n$ generated by
$s_1^1 s_1^2,\dots , s_{n-1}^1s_{n-1}^2$. 

The trivial diagonal idempotent 
$ e_n^{12}= (n!)^{-1} \sum \sigma^{12}$ is such that, for any
$\sigma\in\mfS_n$, one has 
$$ e_n^{12} \sigma^1 = e_n^{12} (\sigma^2)^{-1}
  \quad \& \quad \sigma^1 e_n^{12} = (\sigma^2)^{-1} e_n^{12} \, .$$
Therefore, for any pair $(t,u)$ of tableaux of the same shape,
one has 
$$ e_{t,t}^1\, e_n^{12}\,  e_{u,u}^1  = c\,  e_{t,u}^1\, e_{t,u}^2 \, $$
with some non-zero constant $c$. 

Checking this constant \cite{SgPfaff} furnishes the following
expression of $e_n^{12}$, 
noting as usual $f_\l$ the number of standard tableaux in $\Tab(\l)$~:
\begin{equation}
 e_n^{12}= \sum_{\l:\, |\l|=n} \frac{1}{f_\l} 
    \sum_{t,u \in \Tab(\l)}  e_{t,u}^1 e_{t,u}^2 \, .
\end{equation}

\medskip
Let us now take $\mfS_n^1=\mfS_n$, the second copy $\mfS_n^2$
being the symmetric group $\mfS_n^x$ acting by permutation
on $x_1,\dots,x_n$. Denote by 
$\widetilde \sigma$ the elements of the diagonal group
$\widetilde \mfS_n$. Notice that 
$\Omega_\l = (\sum_{\sigma\in \mfS_n}  \sigma \sigma^x) \, Y_\tc$.

The element $\Omega_\l$ is by construction invariant 
under left multiplication
by all $\widetilde s_i$.  Therefore one has
\begin{equation}   \label{Sum_e_e1}
 \Omega_\l = \widetilde{e_n}\, \Omega_\l
  = \sum_{t,u\in \Tab(\l)   } e_{t,u} e_{t,u}^x \, \Omega_\l 
  =  \sum_{t,u \in \Tab(\l)} e_{t,u}\, (e_{t,u}^x Y_\tc)\, .
\end{equation}  
    
However, Specht polynomials in $V_\l$ are annihilated by all 
$e_\nu^x$, $\nu\neq \mu= \l^{\sim}$. Consequently, the sum
 (\ref{Sum_e_e1}) restricts to 
\begin{equation}   \label{Sum_e_e2}
  \Omega_\l = \widetilde e_\mu \Omega_\l 
 = \sum_{t,u \in \Tab(\l)}  e_{t,u}\, ( e_{t,u}^x\, Y_\tc)  \, .
\end{equation}

On the other hand,
$\Omega_\l s_i =- \Omega_\l$ for all $s_i$ belonging to
the Young subgroup $\mfS_\l$,
because for those $s_i$ one has 
$\Delta(\sigma s_i)= -\Delta(\sigma s_i)$.
Therefore  $\Omega_\l$ is alternating under the
Young subgroup $\mfS_\l$ acting on its right
and in the sum (\ref{Sum_e_e2}),
 the tableaux $u$ must coincide with the first tableau 
$\tc$.  The sum becomes
\begin{equation}   \label{Sum_e_e3} 
\Omega_\l = \sum_{t\in \Tab(\l)} e_{t,\tc} e_{t,\tc}^x\, \Omega_\l
= \frac{n!}{f_\l}  \sum_{t\in \Tab(\l)}  e_{t,\tc} \, 
                \left(e_{t,\tc}^x Y_{\tc}\right)
=  \sum_{t\in\Tab(\l)} \frac{n!}{f_\l}  e_{t,\tc}  Y_t  \, .
\end{equation}  

In final, one has obtained the following expression.

\begin{theorem}
Given a partition $\l$ of $n$, then the element $\Omega_\l$ decomposes as
\begin{equation}  \label{Sum_e_e4}
  \frac{f_\l}{n!}\,  \Omega_\l= \sum_{t\in\Tab(\l)} Y_t\, e_{t,\tc} 
  = \left( \sum_{t\in\Tab(\l)} \zeta^x(t,\tc)
    \zeta(t,\tc) \right) Y_\tc\, e_{\tc,\tc}    \, .
\end{equation}
This is moreover a quasi-idempotent :
$$ \Omega_\l \Omega_\l = \frac{n!}{f_\l} Y_\tc\, \Omega_\l\, .$$
\end{theorem}

For example, for $\l=[2,2]$, one has
$$   \Omega_{2,2} = \frac{4!}{2} 
  \left(  \frac{s_2^x- \frac{1}{2}}{\sqrt{1 - \frac{1}{4}}}
    \frac{s_2 - \frac{1}{2}}{\sqrt{1 - \frac{1}{4}}}
    +1 \right)\, (x_1\moins x_2)(x_3\moins x_4)\, e_{\tc,\tc} \, .$$

Polynomials offer the the possibility of using interpolation methods.
For example, let us consider the specialization
$\x=[0^{\mu_1}, 1^{\mu_2},\dots]$, i.e. $x_1=0=\dots=x_{\mu_1}$,
 $x_{\mu_1+1}=1=\dots=x_{\mu_1+\mu_2},\dots $
The Specht polynomials corresponding to the standard tableaux
of shape $\mu=\l^{\sim}$, which are product of Vandermonde functions,
  vanish uniformly for this specialization,
except in the case of the last tableau $\tr$. Indeed,
there is only one way of distributing the integers $0^{\mu_1}, 1^{\mu_2},\dots$
to avoid having two equal entries in the same column of
the specialized tableau.
For example, 
for $\l=[3,2]$, interpreting each tableau as a product of Vandermonde
functions corresponding to its columns, one has the following
specializations~:

$
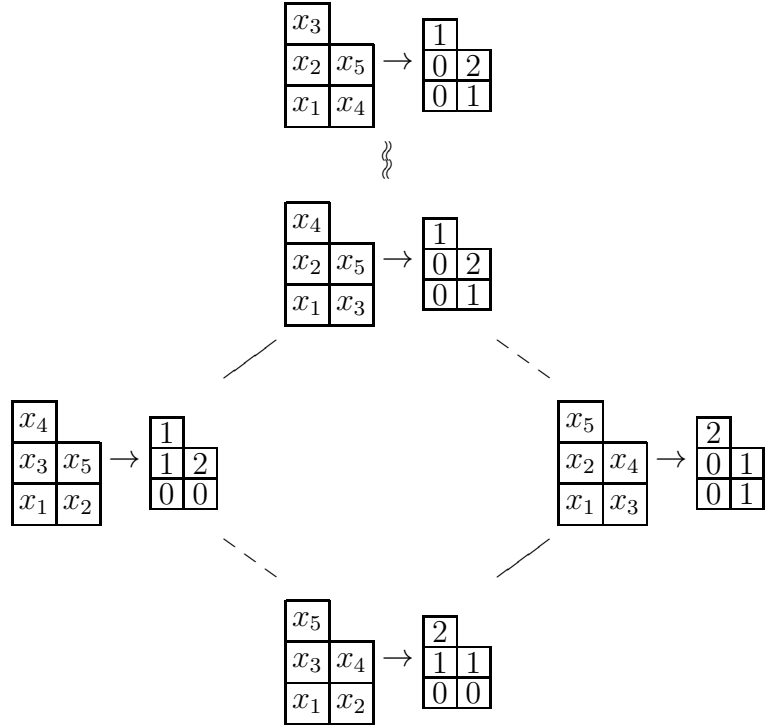
\xymatrix@R=0.5cm@C=0.5cm{
 &  &  & *{\GrTeXBox{
 \young{x_3\cr x_2 &x_5\cr x_1 & x_4\cr}  \to
 \smallyoung{1\cr 0&2\cr 0 &1\cr}
}}\arx5[d]& \\
 &  &  & *{\GrTeXBox{
  \young{x_4\cr x_2 &x_5\cr x_1 & x_3\cr}  \to
 \smallyoung{1\cr 0&2\cr 0 &1\cr}
}}\arx1[ld]\arx3[rd]& \\
 &  & *{\GrTeXBox{
  \young{x_4\cr x_3 &x_5\cr x_1 & x_2\cr}  \to
 \smallyoung{1\cr 1&2\cr 0 &0\cr}
}}\arx3[rd]&  & *{\GrTeXBox{
 \young{x_5\cr x_2 &x_4\cr x_1 & x_3\cr}  \to
 \smallyoung{2\cr 0&1\cr 0 &1\cr}
}}\arx1[ld]& \\
 &  &  & *{\GrTeXBox{
 \young{x_5\cr x_3 &x_4\cr x_1 & x_2\cr}  \to
 \smallyoung{2\cr 1&1\cr 0 &0\cr} 
 }}& \\
}
$
  
Because of the triangularity of the Young basis with respect 
to the Specht basis, one obtains the following proposition.

\begin{proposition}
Let $\l$ be a partition of $n$, $\mu=\l^{\sim}$. 
Then the specialization $\x=[0^{\mu_1}, 1^{\mu_2},\dots]$
of $\Omega_\l$ is equal to 
$$ g_\l := (\l_1\moins 1)! (\l_2\moins 1)! \cdots
   \frac{n!}{\l} e_{\tr,\tc}  \, . $$   

\end{proposition}

Notice that $e_{\tr,\tc}e_{\tr,\tc}=0$, 
the specialization of $\Omega_\l$ does not furnish a quasi-idempotent
because $Y_\tc$ specializes to $0$.
On the other hand, $\sigma_{\tc,\tr} e_{\tr,\tc}$ is such that
$$  \left(\sigma_{\tc,\tr} e_{\tr,\tc} \right)
  \left(\sigma_{\tc,\tr} e_{\tr,\tc} \right)
=   \sigma_{\tc,\tr} e_{\tr,\tc} \, . $$ 

Continuing with the same example, the element $\frac{1}{\s 2} \Omega_{32}$
specializes into $12 e_{\tr,\tc}$ and can be factorized into
$$  g_{3,2} =
s_2s_4s_3 \bigl( 1+s_1s_2s_3 + s_2s_4s_3 + s_2s_3s_4s_1s_2s_3\bigr)
  \, \sum_{\sigma \in\mfS_{3,2}} (\moins 1)^{\ell(\sigma)} \sigma \, .$$
In that case, $\sigma_{\tr,\tc} =s_3s_2s_4$ and one can check
directly that  $s_3 s_2s_4 g_{3,2}$ is a quasi-idempotent.

\bigskip
The coefficients in the expansion of a generic idempotent are non null,
and it is not easy to manipulate sums $\sum c_\sigma\,  \sigma$
of $n!$ terms.  However, one can transform formally such a sum
into a polynomial $\sum c_\sigma y^{\c(\sigma)}$, using the
\emph{code} $\c(\sigma)$ of a permutation defined by
$$ \c(\sigma) =[c_1,\dots, c_n], \quad ,\quad 
 c_i = \#\{ j\ge i,\, \sigma_j \le \sigma_i \} \, .$$

We have exhibited in \cite{YoungAsPol}, for each partition $\l$,
an element in the isotypic component $e_\l \C[\mfS_n] e_\l$,
which, as a polynomial  in $\y$,  factorizes into simple
factors.  In the notation of the present article,
this distinguished element is $\omega g_\l \omega$,
with $\omega=[n,\dots,1]$. 

For example, $[5,4,3,2,1]\, g_{3,2} \, [5,4,3,2,1]$ 
is sent onto
$$ y_1y_2 (1+y_1) (1+y_2)(y_1^2-y_2) (1-y_3+y_3^2) (y_4-1) $$
which represents a sum of $2^4\times 3$ permutations with
coefficients equal to $\pm 1$. A more impressive example 
is the image of $[7,6,5,4,3,2,1]\, g_{4,2,1} \, [7,6,5,4,3,2,1]$,
which is
$$ y^{4220000} (1\plus y_1\plus y_1^2) (1\plus y_2) (1\plus y_3)
(y_2^2 \moins y_3) (1\moins y_4\plus y_4^2 \moins y_4^3)
 (1\moins y_5\plus y_5^2) (1\moins y_6) \, , $$
which expands into a sum of
$2^4\times 3\times 4 =576$ permutations.
We refer to \cite{YoungAsPol} for the precise rule for writing
the factors (we have here an extra monomial factor, because
we are not using the same normalization as in \cite{YoungAsPol}).

\medskip
\begin{center}
\small CNRS, IGM, Universit\'e de Paris-Est\\
\small 77454 Marne-la-Vall\'ee CEDEX 2\\
\small Alain.Lascoux@univ-mlv.fr\\
 \small {\tt http://phalanstere.univ-mlv.fr/$\sim$al}
\end{center}

\end{document}